\documentclass[a4paper,11pt]{article}

\usepackage{amscd}
\usepackage{amsmath}
\usepackage{amssymb}
\usepackage{amsthm}

\theoremstyle{plain}
\newtheorem{theorem}{Theorem}

\theoremstyle{definition}

\newdimen\argwidth
\def\db[#1\db]{%
 \setbox0=\hbox{$#1$}\argwidth=\wd0
 \setbox0=\hbox{$\left[\box0\right]$}
   \advance\argwidth by -\wd0
 \left[\kern.3\argwidth\box0 \kern.3\argwidth\right]}

\newcommand{\diag}{\ensuremath{\mathrm{diag}}}

\newcommand{\Hom}{\mathop{\mathrm{Hom}}\nolimits}

\newcommand{\Spec}{\mathop{\mathrm{Spec}}\nolimits}
\newcommand{\GL}{GL}

\newcommand{\bCx}{\bC^{\times}}

\newcommand{\bC}{\ensuremath{\mathbb{C}}}

\newcommand{\bN}{\ensuremath{\mathbb{N}}}

\newcommand{\bZ}{\ensuremath{\mathbb{Z}}}

\newcommand{\scI}{\ensuremath{\mathcal{I}}}

\newcommand{\scO}{\ensuremath{\mathcal{O}}}

\newcommand{\qgr}{\operatorname{qgr}}
\newcommand{\tor}{\operatorname{tor}}
\newcommand{\gr}{\operatorname{gr}}
\newcommand{\grproj}{\operatorname{grproj}}
\newcommand{\Dbsing}{\ensuremath{D^{\mathrm{gr}}_{\mathrm{Sg}}}}
\newcommand{\Nat}{\ensuremath{\mathrm{Nat}}}
\newcommand{\module}{\operatorname{mod}}

\title{Triangulated categories
of Gorenstein cyclic quotient singularities}
\author{Kazushi Ueda}
\date{}
\begin{document}

\maketitle

\begin{abstract}
We prove an equivalence of triangulated categories
between Orlov's triangulated category of singularities
for a Gorenstein cyclic quotient singularity
and the derived category of representations
of a quiver with relations
which is obtained from the McKay quiver
by removing one vertex and half of the arrows.
\end{abstract}

Fix an integer $n$
greater than one.
For a finite subgroup $G$ of $\GL_n(\bC)$,
let $\{ \rho_i \}_{i=0}^{N-1}$ be the set of
irreducible representations of $G$.
Let further
$\rho_{\Nat}$ be
the natural $n$-dimensional representation
of $G$
given by the inclusion.
For $k, l = 0, \dots, N-1$,
define the natural numbers
$a_{kl}$
by the decomposition
$$
 \rho_l \otimes \rho_{\Nat} = \bigoplus_k \rho_k^{\oplus a_{kl}}
$$
of the tensor product
into the direct sum of irreducible representations.
The {\em McKay quiver} of $G$
is the quiver whose set of vertices is
$\{ \rho_i \}_{i=0}^{N-1}$
and the number of whose arrows
from the vertex $\rho_k$ to the vertex $\rho_l$
is $a_{kl}$
\cite{McKay_GSFG}.

Now assume that
$G$ is a cyclic group
generated by an element of the form
$
 g = \diag(\zeta^{a_1}, \cdots, \zeta^{a_n}),
$
where
$a_1, \ldots, a_n$ are positive integers
such that $\gcd(a_1, \ldots, a_n) = 1$ and
$\zeta = \exp[2 \pi \sqrt{-1}  / (a_1 + \cdots + a_n)]$.
Let $R = \bC[x_1, \ldots, x_n]$ be the polynomial ring
in $n$ variables
equipped with the $\bZ$-grading
given by $\deg x_i = a_i$, $i = 1, \ldots, n$.
Define another $\bZ$-graded ring
$A(a_1, \dots, a_n) = \bigoplus_{k \geq 0} A_k$
by
\begin{equation} \label{eq:A}
 A_k = R_{k (a_1 + \cdots + a_n)}.
\end{equation}
Then $A(a_1, \dots, a_n)$ is the invariant ring of $R$
by the action of $G$
so that $\bC^n / G = \Spec A(a_1, \dots, a_n)$.
%

In this case,
the corresponding McKay quiver has
$N = a_1 + \cdots + a_n$ vertices
$\{ \rho_{k} \}_{k=0}^{N-1}$
and $n N$ arrows
$\{ x_{i, k} \}_
{\tiny 
\begin{array}{l}
 i=1, \dots, n \\
 k=0, \dots, N-1
\end{array}
},
$
where $\rho_k$
is the one-dimensional representation of $G$
sending $g \in G$ to
$
 \exp[- 2 k \pi \sqrt{-1} / (a_1 + \cdots + a_n)] \in \bCx,
$
and $x_{i, k}$ is an arrow
from $\rho_k$ to $\rho_{k + a_i}$.

Next we introduce another quiver $Q(a_1, \dots, a_n)$
obtained
by removing the vertex $\rho_0$
and half of the arrows
from the McKay quiver;
the set of vertices of $Q(a_1, \dots, a_n)$
is $\{ \rho_k \}_{k=1}^N$,
and an arrow of the McKay quiver
from $\rho_k$ to $\rho_l$
is an arrow of $Q_g$
if $0 < k < l$.

A {\em quiver with relations}
is a pair $\Gamma = (Q, \scI)$
of a quiver $Q$ and
a two-sided ideal $\scI$
of its path algebra $\bC Q$.
We equip $Q(a_1, \dots, a_n)$ with the relations $\scI(a_1, \dots, a_n)$
generated by
$
 x_{j, k + a_i} x_{i, k} - x_{i, k + a_j} x_{j, k}
$
for $1 \leq i < j \leq n$ and $k = 1, \ldots, N-a_i-a_j-1$,
and
put $\Gamma(a_1, \dots, a_n) = (Q(a_1, \dots, a_n), \scI(a_1, \dots, a_n))$.
The main theorem is:

\begin{theorem}
For a sequence $a_1, \dots, a_n$
of positive integers
such that $\gcd(a_1, \dots, a_n) = 1$,
there exists an equivalence
$$
 \Dbsing(A(a_1, \dots, a_n))
  \cong D^b \module \Gamma(a_1, \dots, a_n)
$$
of triangulated categories.
\end{theorem}

Here, $D^b \module \Gamma(a_1, \dots, a_n)$
is the bounded derived category
of finite-dimensional right modules
over the path algebra
$\bC \Gamma(a_1, \dots, a_n)
   = \bC Q(a_1, \dots, a_n) / \scI(a_1, \dots, a_n)$
with relations.

$\Dbsing(A(a_1, \dots, a_n))$ is the
{\em triangulated category of singularities},
defined by Orlov \cite{Orlov_DCCSTCS}
as the quotient category
\begin{equation} \label{eq:Dbsing}
 \Dbsing(A(a_1, \dots, a_n))
  = D^b \gr A(a_1, \dots, a_n) / D^b \grproj A(a_1, \dots, a_n)
\end{equation}
of the bounded derived category
$D^b \gr A(a_1, \dots, a_n)$
of finitely-generated $\bZ$-graded
$A(a_1, \dots, a_n)$-modules
by its full triangulated subcategory
$D^b \grproj A(a_1, \dots, a_n)$
consisting of perfect complexes.
The $n=2$ case
in the above theorem
is due to Takahashi \cite{Takahashi_MF}
(see also
Kajiura, Saito, and Takahashi
\cite{Kajiura-Saito-Takahashi_MF}).

The proof goes as follows:
Let
\begin{equation*}
 \qgr R := \gr R / \tor R
\end{equation*}
be the quotient category
of the abelian category $\gr R$
of finitely-generated
$\bZ$-graded $R$-modules
by its full subcategory $\tor R$ 
consisting of torsion modules,
and $\pi: \gr R \to \qgr R$ be
the natural projection functor.
For $M \in \gr R$ and $l \in \bZ$,
$M(l)$ denotes the graded $R$-module shifted by $l$;
$M(l)_k = M_{l+k}$.
Define a shift operator
$
 s: \qgr R \to \qgr R
$
by
$s(\pi M) = \pi M(a_1+\dots+a_n)$
and put
$
 \scO = \pi R.
$
Then one has
$
 A(a_1, \dots, a_n) = \bigoplus_{k=0}^{\infty} \Hom(\scO, s^k(\scO)).
$
Since $\gcd(a_1, \dots, a_n)=1$,
the graded $R$-module $R(l)$ for any $l \in \bZ$
is generated up to torsion modules
by the subset
$\bigcup_{j \in \bN} R(l)_{j(a_1 + \dots + a_n)}$
consisting of elements whose degrees are
positive multiples of $a_1 + \dots + a_n$.
Hence $s$ is ample and
one has
$$
 \qgr R \cong \qgr A(a_1, \dots, a_n)
$$
by Artin and Zhang \cite[Theorem 4.5]{Artin-Zhang_NPS}.

Since $s^{-1}(\scO)$ is the dualizing sheaf,
$A(a_1, \dots, a_n)$ is Gorenstein with Gorenstein parameter $1$
(cf. \cite[Lemma 2.11]{Orlov_DCCSTCS}).
Therefore one has a semiorthogonal decomposition
$$
 D^b \qgr A(a_1, \dots, a_n)
  = \langle \scO, \Dbsing(A(a_1, \dots, a_n)) \rangle
$$
by
Orlov \cite[Theorem 2.5(i)]{Orlov_DCCSTCS}.
Here,
$D^b \qgr A(a_1, \dots, a_n)$ denotes the bounded derived category
of the abelian category $\qgr A(a_1, \dots, a_n)$.
On the other hand,
$D^b (\qgr R)$ has a full strong exceptional collection
$(\scO, \scO(1), \dots, \scO(a_1 + a_2 + \dots + a_n - 1))$
(see 
e.g.
\cite[Theorem 2.12]{Auroux-Katzarkov-Orlov_WPP}).
Hence
$\Dbsing(A(a_1, \dots, a_n))$ is equivalent
to the full triangulated subcategory
of $D^b \qgr A(a_1, \dots, a_n)$
generated by the exceptional collection
$(\scO(1), \dots, \scO(a_1 + \dots + a_n - 1))$.
By Bondal \cite[Theorem 6.2]{Bondal_RAACS},
this subcategory is isomorphic
to the derived category
of finite-dimensional right modules
over the total morphism algebra
$$
 \bigoplus_{i,j=1}^{a_1 + \dots + a_n - 1}
  \Hom(\scO(i), \scO(j))
$$
of this collection,
which is isomorphic to $\bC \Gamma(a_1, \dots, a_n)$.

{\bf Acknowledgment:}
We thank Kentaro Nagao and the anonymous referee
for pointing out errors
in the earlier version.
This research is supported by
the 21st Century COE Program of Osaka University
and Grant-in-Aid for Young Scientists (No.18840029).

\bibliographystyle{plain}
\newcommand{\noop}[1]{}\def\cprime{$'$} \def\cprime{$'$}
  \def\cydot{\leavevmode\raise.4ex\hbox{.}} \def\cprime{$'$} \def\cprime{$'$}
  \def\cprime{$'$}

\noindent
Kazushi Ueda

Department of Mathematics,
Graduate School of Science,
Osaka University,
Machikaneyama 1-1,
Toyonaka,
Osaka,
560-0043,
Japan.

{\em e-mail address}\ : \  kazushi@math.sci.osaka-u.ac.jp

\end{document}